\documentclass[a4paper,11pt]{amsart}
\usepackage{amscd,amsmath,amssymb,amsxtra,artov2}

\usepackage[xdvi]{graphicx}
\usepackage[]{hyperref}
\def\2{\partial}

\def\vare{\varepsilon}

\def\si{\sigma}

\def\lan{\langle}
\def\ran{\rangle}
\title{Configurations of Skew Lines}

\email{oleg@math.uu.se}
\author{Julia Viro}
\email{julia@math.uu.se}

\author{Oleg Viro}
\address{Department of Mathematics, Uppsala University,
Box 480, S-751 06 Uppsala, Sweden
}

\begin{document} 
\begin{abstract} This article is a survey of 
results on projective configurations of subspaces in general 
position.  It is written in the form of introduction to the 
subject, with much of the material accessible to advanced high school 
students.  However, in the part of the survey concerning configurations 
of lines in general position in three-dimensional space we give a 
detailed exposition.  \end{abstract}

\maketitle

The first version of this paper was written as an elementary introductory text for 
high-school students. It was published \cite{1} in the the journal 
``Kvant'', the third issue of 1988, but in a shortened form. Then we 
expanded the article in order to encompass or at least mention some 
related questions. However we decided to keep the style of 
\cite{1}, in the hope that it would also be appreciated by a 
professional mathematician. We apologized to a reader, who would find 
the style irritating,  and we mentioned that the material in the first 
two-thirds of the article (through the section on ``sets of five 
lines'') was announced in the note \cite3, while the final third of 
the article is written in a more traditional style. 

The expanded version \cite{2} of \cite{1} was published in the first 
volume  of a Russian journal Algebra i Analiz opening a new section 
``Light reading for the professional''. English version of the paper 
became available in a translation made by N. Koblitz and published by 
American Mathematical Society in the first volume of Leningrad 
Mathematical Journal.

Unfortunately, the first volumes, even of the first rate journals, are 
not distributed as well as they deserve. Leningrad Mathematical Journal
is an excellent journal, but we would like to bring our paper to more 
readers. During 16 years which passed since the time of writing \cite{2} some 
questioned posed in \cite{2} were solved and we decided to refresh the 
text and make it available to new readers. Partly we fix here some 
of defects of \cite{2}: wrong pictures and few terminological 
inaccuracies of translation.

Furthermore, some of the results appeared in other papers without appropriate references. The tradition of wrong priority references was established in papers
by R.\ Penne and H.\ Crapo \cite{P1}, \cite{CP}. We expanded the bibliography 
for the sake of completeness, but arranged it in the chronological order, 
to facilitate orientation in the history of the subject.

\section*{Can skew lines be interlaced?}
\par
The article in ``Kvant'' was titled {\it Interlacing of Skew Lines\/}. This
title sounds a little strange, doesn't it? The word ``interlacing'' suggests
something flexible, not straight lines! To be sure, the title refers not to be
the process of interlacing, but rather to the result. But is it possible to
weave together skew lines which are situated in some clever way with respect to
one another? At first glance this may seem not to be possible. Yet where do we
get this impression? In daily life we never come across anything that really
resembles a straight line. What bothers us is not that there is not such thing
as an infinitely thin object---we are prepared to neglect the thickness---but
rather that there is no such thing as an infinitely long object. Even light
rays---which are models of linearity---become
 scattered and dispersed, and cannot be detected at a large distance. In
practice one deals only with line segments.
\par
Any set of disjoint line segments can be moved around to any other relative
location in such a way that they remain disjoint. This we can see from
experience, and it is also not hard to prove. We depict straight lines using
line segments, and so it seems to us that straight lines cannot be woven
together. But is that really the case?
\par
First of all, let us give a more precise statement of the questions which
concern us. The first question is: Can a set of disjoint lines be rearranged?
But what do we mean by the term ``rearrange''? Here we shall not be concerned
with the angles or distances between the lines. We shall consider the relative
position of the lines to be unchanged if we move them in such a way that they
never touch. But if one set of lines cannot be obtained from another set by
such a movement, then we shall say that the two sets of lines are arranged
differently.
\par
The simplest lines for us to visualize are parallel lines. Clearly, any two
sets of parallel lines with the same number of lines in each set have the same
arrangement. In fact, if we consider the lines of one set to be ``frozen'' in
place and then rotate the entire space, we can make them parallel to the lines
of the other set; then, moving the lines of the first set one by one in such a
way that they remain parallel and do not bump into one another, we can easily
make them coincide with the lines of the second set.
\par
We now consider arbitrary sets of lines. Can an arbitrary set of lines be moved
(``combed'') into a set of parallel lines? This question has a simple and
unexpected answer, which is hard to arrive at by considering concrete sets of
lines. If you take a specific set of lines and study it for a while, you can
probably find a way to make all of the lines parallel. But this does not give
an answer to the question in full generality, because you undoubtedly made use
of some specific features of your set of lines. Can one treat all possible sets
of lines at once? It turns out that one can, and this is how. Let us take an
arbitrary set of disjoint lines. We choose two parallel planes which are not
parallel to any of the lines in our set. We fix the points of intersection of
the first plane with the lines, fastening the lines at those points. We also fix
the intersection of the lines with the second plane, but only as a point on that
plane, which we allow to slide along the lines. In other words, we drill small
holes in the second plane where it intersects with the lines. We then move the
second plane away from the first one in the direction perpendicular to both
planes. The lines are pulled through the little holes, and the angles which
they form with the planes increase. If we move the second plane to infinity in a
finite amount of time, then these angles all reach $90^\circ$, i.e., the lines
become parallel to one another. This ``combing'' of our set of lines can be
described as follows in a language which is more customary for geometry: we
expand the space away from the first plane in a direction perpendicular to it,
where the expansion factor increases rapidly to infinity in a finite length of
time. Here the straight lines rotate around their points of intersection with
the plane, and in the limit they become perpendicular to the plane.
\par
Thus, one cannot have interlaced disjoint lines: all sets of disjoint lines
have the same arrangement. But our title refers to skew lines, and so sets of
parallel lines are excluded. There is a serious reason for this. Parallel lines
are very close to being intersecting lines: one can move one of two parallel
lines by an arbitrarily small amount so as to make them intersect. This is not
the case for skew lines.
\par
Since we have decided not to allow parallel lines, we must reexamine the
question of which sets of lines have the same arrangement and which do not. We
shall say that the arrangement of a set of lines remains the same if it is
moved in such a way that the lines are always skew, never parallel. In what
follows we will often be considering such movements of lines, and so it is
useful to have a special word to refer to them. We shall use the word {\it
isotopy\/} to denote such a movement of lines. If one set of lines cannot be
obtained from another by means of an isotopy, then we say that the two sets
have different arrangements. We shall also say that such sets of lines are
nonisotopic.
\par
The amount of difficulty in determining whether two sets of lines are isotopic
depends most of all on the number of lines in the sets. In general, the more
lines, the more clever one must be to find an isotopy which transforms one set
into the other. We first treat the simplest case of the isotopy problem.
 
\section*{ Two Lines}
We take any two pairs of skew lines, and try to decide whether they are
isotopic. In this case it is perhaps too pretentious to use the word
``problem'', because it is completely obvious that we have an isotopy.
Nevertheless, we shall make a detailed examination of the proof.
\par
Using a rotation around a line which is perpendicular to both lines in one of
the pairs, we can make the angle between the lines the same in both pairs; in
fact, we can make both angles $90^\circ$. We note that the smallest line
segment joining the two lines in a pair is the segment of the common
perpendicular which is contained between them. We next bring the two lines
closer together (or move them farther apart) along this perpendicular, so that
the segments have the same length for the two pairs; after that we move one
pair so that the segment between the two lines coincides with the segment for
the other pair.
 We use a rotation around this segment to make one of the lines of the first
pair coincide with a line of the second pair (this can be done because all of
the lines are perpendicular to the segment). In the process the second lines of
the pairs also come together. In fact, they both pass through a common
point---an endpoint of the perpendicular segment---and are perpendicular to the
same plane---the plane determined by the perpendicular and the first lines of
the pairs (which now coincide). The proof is complete.
\par
\begin{figure}[htb]
\centerline{\includegraphics{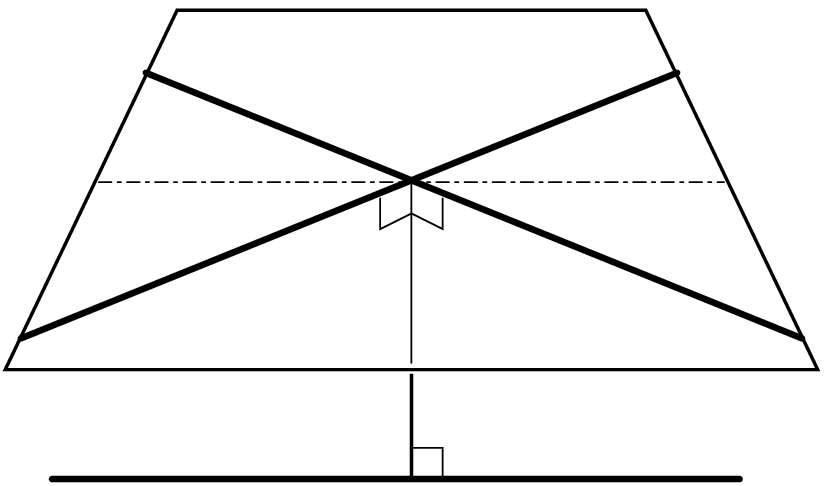}}
\caption{}
\label{f1}
\end{figure}

At the end of the proof, after we made the distances between the two lines the
same for the two pairs, we moved a pair of lines in a rigid manner---without
changing either the distance or the angle between them. The question arises:
Suppose that both the distances and angles between the two lines are the same
for two pairs of skew lines. Is it always possible to find an isotopy between
the two pairs during which the distance and angle remain fixed? The previous
argument shows that this question has an affirmative answer if the angle is
$90^\circ$. However, if the angle is not $90^\circ$, then it may happen that
after the isotopy in the previous paragraph the second lines in the pairs do
not coincide. This unlucky
 case is illustrated in Figure \ref{f1}. The second lines in
the pairs form an angle whose bisector is parallel to the first (skew) lines,
and the plane containing the second lines is perpendicular to the plane
containing the bisector and the first lines. Thus, there was a good reason why
we wanted to make the angles $90^\circ$ in the beginning of the above proof:
for any other choice of the angle, the construction would not give the desired
result. But this was not simply an artifact of our particular construction; it
turns out that any two pairs of skew lines with equal distance and angle which
do not coincide after the above construction cannot be made to coincide using
any isotopy during which the distance and angle remain fixed. This is connected
with a remarkable phenomenon, which we shall encounter often in the sequel. It
merits a more detailed discussion.

\section*{Orientations and Semi-Orientations}

\begin{figure}[htb]
\centerline{\includegraphics{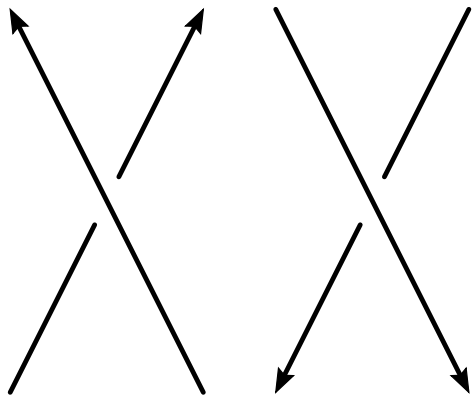}}
\caption{Semi-orientation of a pair of lines.}
\label{f2}
\end{figure}

To orient a set of lines means to give a direction to each line in the set. There
are $2^n$ possible orientations of a set of $n$ lines. A {\it
semi-orientation\/} of a set of lines is a pair of opposite orientations (Figure
\ref{f2}).
\par
\begin{figure}[htb]
\centerline{\includegraphics{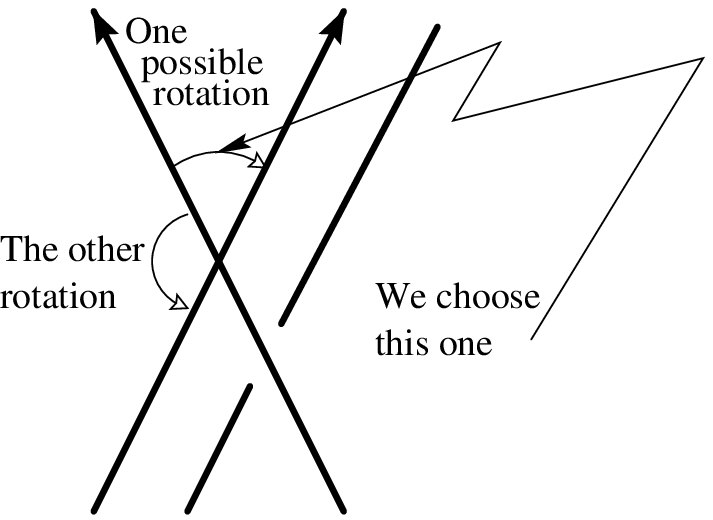}}
\caption{}
\label{f3}
\end{figure}
\begin{figure}[hbt]
\centerline{\includegraphics{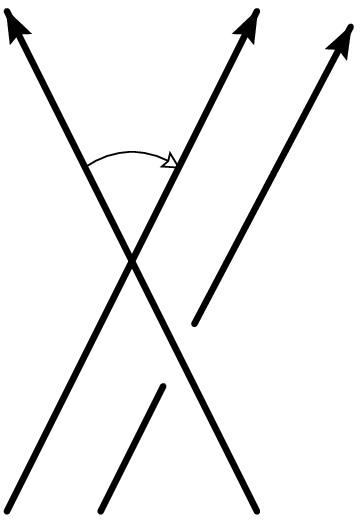}}
\caption{}
\label{f4}
\end{figure}

Any pair of nonperpendicular lines has a canonical semi-orientation which is
determined by the relative position of the two lines. Namely, we choose an
arbitrary orientation of one of the lines, and then we determine the orientation
of the second line by rotating the first line in the most economical way (i.e.,
with the smallest angle of rotation) so as to make it parallel to the second
line (see Figure \ref{f3}). We then give the second line the 
orientation pointing in
the same direction as the (now parallel) first line (Figure \ref{f4}). 
Thus, choosing
an orientation of one of the lines determines an orientation of the pair. If we
choose the opposite orientation of the first line, then we obtain the opposite
orientation of the pair. If we were to use the other line to start with, we
would obtain the same pair of opposite orientations.
These two opposite orientations are what we meant by the canonical
semi-orientation of the pair of nonperpendicular lines.

\par
An isotropy during which the angle between the lines remains fixed takes the
canonical semi-orientation to the canonical semi-orientation. This suggests the
idea of considering another type of isotopy---isotopies of semi-oriented pairs
of skew lines. Here we allow the angle and distance between the lines to change,
but we require that the semi-orientation be preserved. Such an isotopy occupies
an intermediate position between an arbitrary isotopy and an isotopy during
which the distance and angle (where we suppose that the angle is
$\not=90^\circ$) remain fixed. That is, if there is no semi-oriented isotopy
between two semi-oriented pairs of lines, then there is certainly no isotopy
between them which preserves the distance and angle. What can stand in the way
of an isotopy of semi-oriented pairs of lines?
\section*{The Linking Number}

\begin{figure}[htb] 
\centerline{\includegraphics{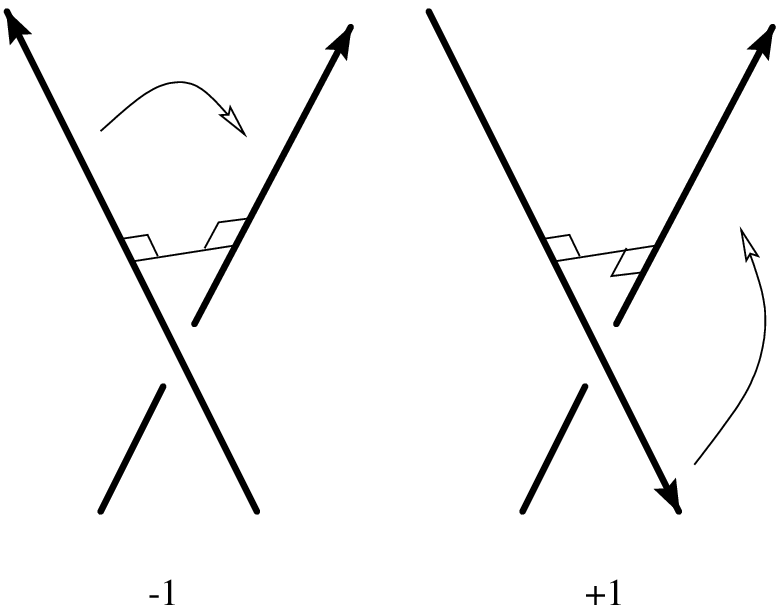}}
\caption{}
\label{f5}
\end{figure}

Any semi-oriented pair of lines has a characteristic which takes the value $+1$
or $-1$. It is called the {\it linking number\/}. This number is
preserved under isotopies, and so if two semi-oriented pairs of lines have
different linking numbers, then they are not isotopic. Here is the
definition of the linking number. The most economical way of aligning an
oriented line with a second oriented line which is skew to it is to place it
alongside a common perpendicular to the two lines and then rotate it by the
smallest angle that brings it to the same direction as the second line. Here
the line rotates either like the right hand around the thumb, or like the left
hand (Figure \ref{f5}). In the first case the linking number is $-1$, 
and in the
second case it is $+1$.

\par

\begin{figure}[htb]
\centerline{\includegraphics{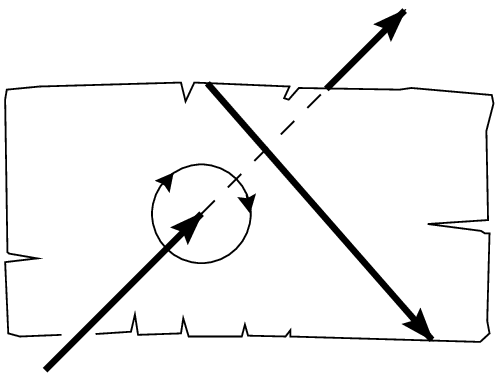}}
\caption{}
\label{f6}
\end{figure}

\begin{figure}[htb]
\centerline{\includegraphics{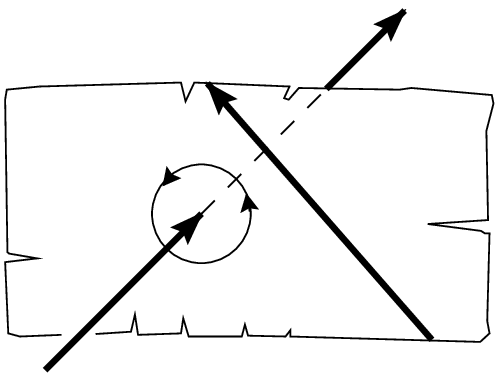}}
\caption{}
\label{f7}
\end{figure}

\begin{figure}[hbt]
\centerline{\includegraphics{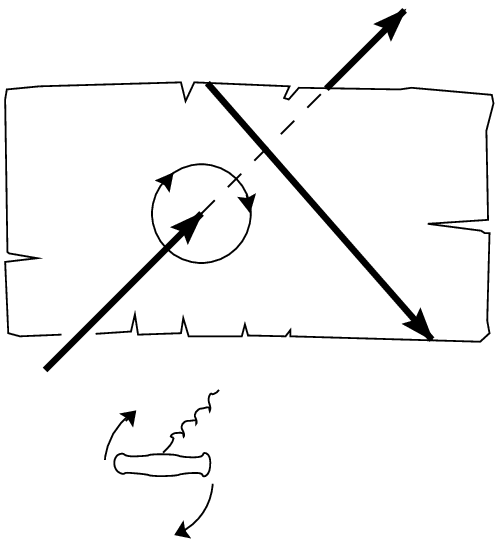}}
\caption{}
\label{f8}
\end{figure}

To help the reader familiar with algebraic topology make the right connection,
we give a second equivalent definition of the linking number of a pair of
oriented skew lines. Through one of the lines we draw a plane which intersects
the other line. We place our right hand so that our thumb rests on the second
line and passes through the plane in the direction determined by the orientation
of the line, while rotating in the direction our fingers point. On the plane we
obtain an oriented circle which is traced by the tips of our fingers. The
orientation of the circle may be the same as the orientation of the first line
(Figure \ref{f6}) or different (Figure \ref{f7}). In the first case 
the linking number
is $+1$, and in the second case it is $-1$. Figure \ref{f8} will 
enable the reader to
see that the two definitions of the linking number are equivalent.

\par
\begin{figure}[hbt]
\centerline{\includegraphics{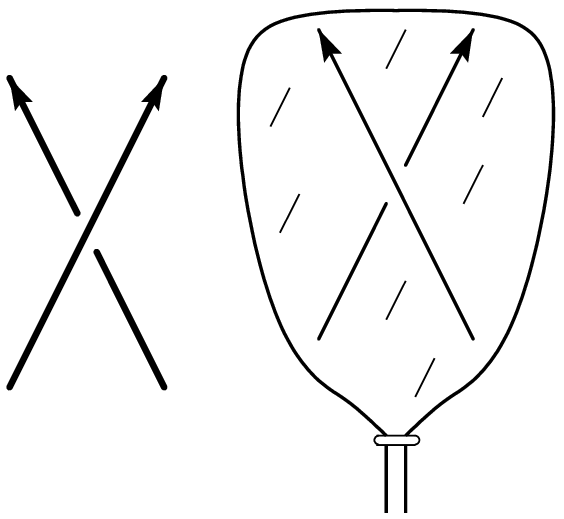}}
\caption{}
\label{f9}
\end{figure}

It is clear that changing the orientation of one of the lines of the pair
changes the linking number. Hence, if the orientation of the pair is
changed to the opposite orientation (i.e., the orientation is reversed on both
lines), then the linking number does not change. In other words, the
linking number is an invariant of a semi-oriented pair: it depends only on
the semi-orientation. If we look at the reflection of our pair of oriented
lines in the mirror (Figure \ref{f9}), the linking number 
changes. \par

\begin{figure}[htb]
\centerline{\includegraphics{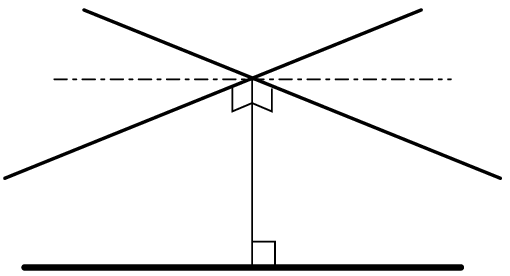}}
\caption{}
\label{f10}
\end{figure}

\begin{figure}[htb]
\centerline{\includegraphics{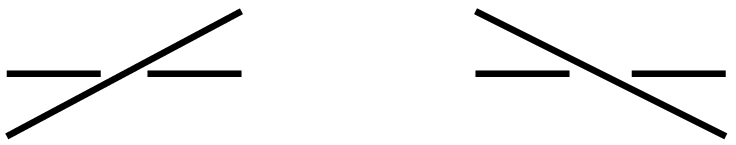}}
\caption{}
\label{f11}
\end{figure}

We now return to the unfortunate situation we encountered when looking for an
isotopy between two pairs of skew lines which preserves the distance and angle
between the lines (see Figure \ref{f10}). At the time we could not 
answer the question
of whether the sets are isotopic (Figure \ref{f11}). Now, however, we 
see that these
pairs (with their canonical semi-orientation) are obtained from one another by
a mirror reflection, and so they have different linking numbers. Thus,
they cannot be connected by an isotopy which preserves the distance and angle
between the lines. But if two pairs have the same distance and angle and also
the same linking number, then they can be connected by such an isotopy.

By the way, it is possible to modify the notion of the angle between two skew
lines in such a way as to incorporate the linking number and thereby make
it unnecessary to work with the linking number separately. The angle
between two lines was defined above so as to be in the interval
$(0^\circ,90^\circ)$.
 define the modified angle between two lines to be the
product of the angle in the earlier sense and the linking number, if the
latter is defined (i.e., if the angle is not $90^\circ$), and to be the angle
in the earlier sense (i.e., $90^\circ$) if the linking number is not
defined. The modified angle is in one of the intervals $(-90^\circ,0^\circ)$,
$(0^\circ,+90^\circ]$. The sign can be determined from the right hand rule,
without saying anything about the linking number.
\par
We have thereby completely analyzed the situation with sets of two skew lines.

\section*{Triples of Lines}
\par
When we studied pairs of lines, an important role was played by the common
perpendicular to the two skew lines. Strictly speaking, we could have avoided
using it; but it seemed to be connected to the lines in such a natural way,
providing a tangible tie between them, that it would have been strange not to
make use of it. Now it would be good to find something equally natural for a
triple of skew lines. There are two objects that are capable of playing this
role. We shall discuss one of them now, and postpone consideration of the
second one. Jumping ahead, suffice it to say that the second object is a
hyperboloid.
\par
These objects will not be associated to every triple of pairwise skew lines. We
will have to disallow triples whose lines lie in three parallel planes. But
notice that such an arrangement is unstable: by nudging one of the lines a
little, we obtain an isotopic triple to which our constructions can be applied.
\par

\begin{figure}
\centerline{\includegraphics{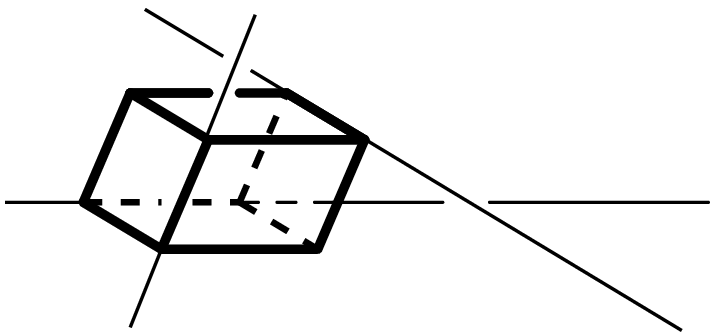}}
\caption{}
\label{f12}
\end{figure}

Thus, we consider an arbitrary triple of pairwise skew lines which do not lie
in three parallel planes. For each line we draw two planes containing the line,
each parallel to one of the other two lines. In this way we obtain six planes,
i.e., three pairs of parallel planes. These planes intersect to form a
parallelepiped. Our lines are the extensions of three of its skew edges (Figure
\ref{f12}). Thus, any three pair-wise skew lines which do not lie in 
three parallel
planes are extensions of the edges of a certain parallelepiped. This
parallelepiped is the first object which we associate to the triple of lines.
What is special about it? In the first place, it is unique. In fact, there is a
unique plane parallel to a given line that contains a second skew line; and if
these lines are the extensions of edges of a parallelepiped, then this plane
contains one of its faces. Consequently, the six planes are uniquely determined
by the original triple of lines; since any parallelepiped whose edges lie on
these lines is bounded by those planes, it is also uniquely determined.

We see that the parallelepiped joins together the lines of the triple just as
nicely as the common perpendicular joined together the lines of the pair. Just
as in the case of the common perpendicular and the semi-oriented pair of
nonperpendicular lines, the original geometry of the configuration naturally
leads to something more, though still something which is connected with it in a
canonical way and so merits our further consideration when we study the
original object.
\vskip .1in

\begin{figure}
\centerline{\includegraphics{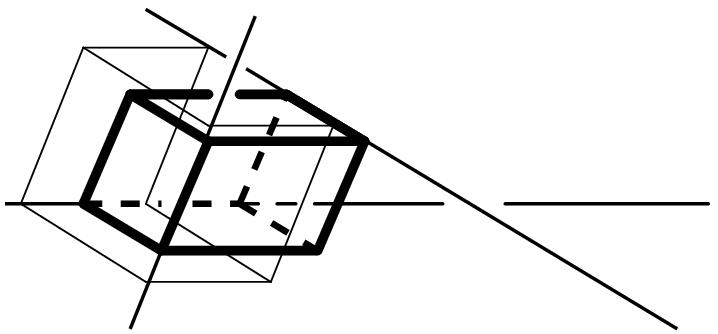}}
\caption{}
\label{f13}
\end{figure}

{\bfit A riddle} In Figure \ref{f13}, despite what was proven 
above, we have drawn two
different parallelepipeds with edges lying on three pairwise skew lines. What is
going on?

\begin{figure}
\centerline{\includegraphics{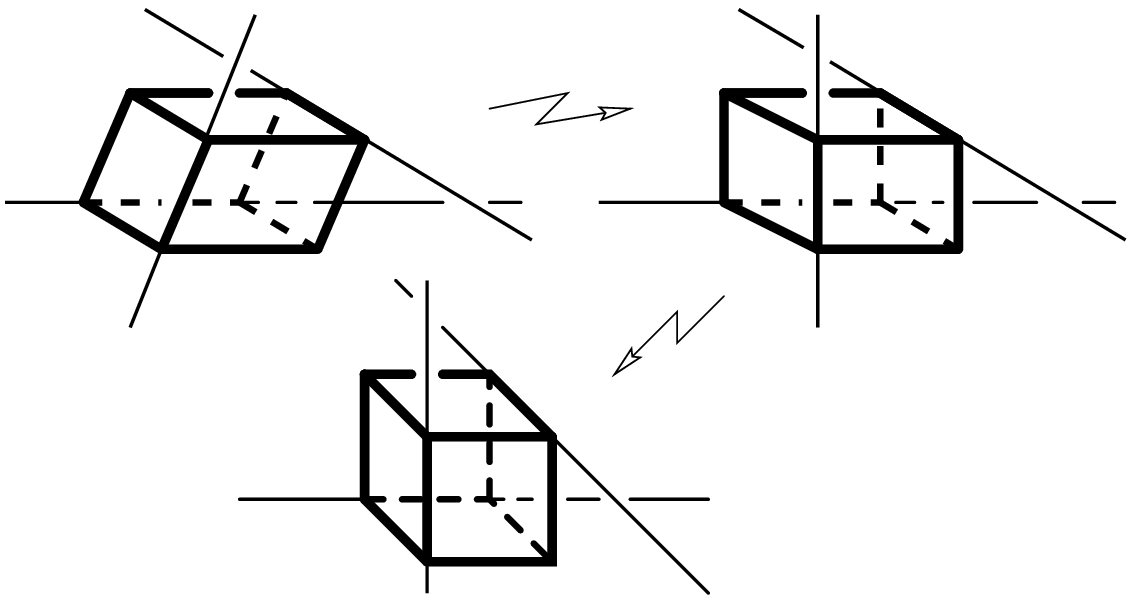}}
\caption{}
\label{f14}
\end{figure}

\begin{figure}
\centerline{\includegraphics{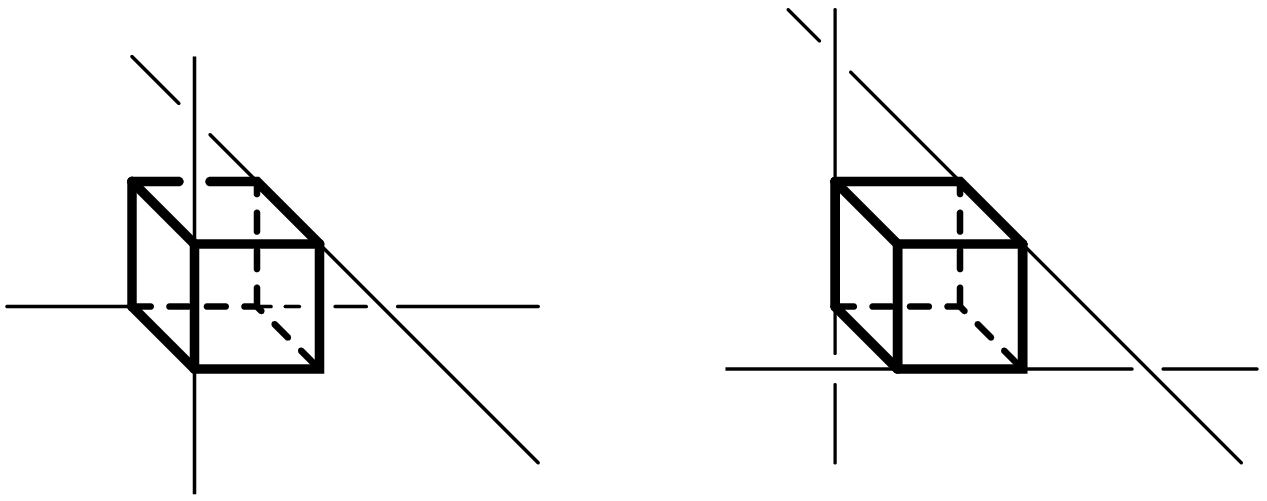}}
\caption{}
\label{f15}
\end{figure}

We now look at the classification of triples up to isotopy. As shown above, we
may suppose that the lines in the triple are extensions of edges of a certain
parallelepiped. A parallelepiped is determined (up to translation) by the
lengths of its edges and the angles between them. Using a continuous
deformation, we can first make all of the angles into right angles (obtaining a
rectangular parallelepiped), and then we can make all of the edges have the
same length, for example, length one (obtaining a cube) (Figure 
\ref{f14}). This
deformation induces an isotopy of the triple of lines which are extensions of
edges of the parallelepiped. In this way we have managed to place the lines of
our triple along pairwise skew edges of a unit cube. This is a remarkable
accomplishment. It means that we now know that there are not very many possible
nonisotopic sets of three skew lines---there are at most the number of triples
of skew edges on a cube, and this number is 8. And even 8 is too many. We can
use a rotation of the cube to take any edge of the cube to any other edge, and
this reduces the number of possible nonisotopic configuration types to two. They
are shown in Figure \ref{f15}.

\par
This success might prompt us to hope that we can similarly find an isotopy
between the two triples of lines in Figure \ref{f15}, and thereby 
prove that all triples of skew lines are isotopic. Try to do this!
\par
You're having trouble? Don't blame yourself---it cannot be done! Just like
pairs of oriented lines, triples of (nonoriented) lines have an invariant, also
called the {\it linking number\/}, which takes the value $+1$ or $-1$, is
preserved under isotopies, and changes when one takes a mirror reflection of
the triple of lines. Here is its definition. Suppose we have a set of three
pairwise skew lines. We orient the tree lines in an arbitrary way. Then each
pair of lines in the triple has a linking number (equal to $\pm1$). If we
multiply all of the linking
 numbers, we obtain a number (also $\pm1$),
which is what we call the linking
 number of the original triple of lines.
This number does not depend on the orientation of the lines, since if we
reverse the orientation of any line, the effect is to change the linking
numbers of two of the pairs, and this does not change the product. The
fact that the linking number of a triple is preserved under isotopy and
changes under mirror reflection follows from the corresponding properties of
the linking numbers of pairs of oriented lines. Since the triples of lines
in Figure \ref{f15} are the mirror images of one another, they have 
different linking
numbers, and hence they are not isotopic to one another.
\par
Since any triple of pairwise skew lines is isotopic to one of the two triples
in Figure \ref{f15}, it follows that two triples of lines are isotopic 
if and only if they have the same linking number.
\par
Thus, as soon as we reach three lines we find that there are different possible
arrangements of triples of skew lines. This provides a justification for the
title of the paper, and for our subsequent use of the word {\it interlacing\/}
for a set of pairwise skew lines.

{\bfit Problem.\/} It is natural to expect that the linking number of 
a triple
of nonoriented lines is equal to the linking number of some pair of
semi-oriented lines which can be constructed from the triple in a canonical
way. This is in fact the case, except that rather than one such semi-oriented
pair there are three. Prove that for any triple of skew lines there is a unique
semi-orientation such that the linking numbers of all three pairs of lines
in the triple are equal. Obviously, this value is also equal to the linking
number of the triple.

\section*{Amphicheiral and Nonamphicheiral Sets}

Note that a triple of skew lines is never isotopic to its mirror image,
while a pair of lines is isotopic to its mirror image. In general, we say that
a set of pairwise skew lines is {\it amphicheiral\/}  if it 
is isotopic
to its mirror image; otherwise we say it is {\it nonamphicheiral.\/}  Thus, a
triple is always nonamphicheiral, and a pair is amphicheiral. The 
following questions arise:
\par
1) Are there other values of $p$ such that any interlacing of $p$ lines is
nonamphicheiral?
\par
2) Are there other values of $p$ such that any interlacing of $p$ lines is
amphicheiral?
\par
3) For what $p$ does there exist a nonamphicheiral interlacing of $p$ 
lines? \par
4) For what $p$ does there exist an amphicheiral interlacing of $p$ 
lines? \par
Although this does not take us very far in the direction of an answer to our
original question (of describing the set of interlacings of $p$ lines up to
isotopy), it is worthwhile to take up these four questions. They are rough and
somewhat superficial questions, but at the same time they have a more
qualitative character. Because of this roughness and superficiality we can be
confident of early success, and the result will undoubtedly be useful in our
classification.
\par
We do not yet have at our disposal very many tools for proving that 
a set is nonamphicheiral. But we do know that every triple is nonamphicheiral, and this is already a
lot. After all, any set of more than three lines contains triples. Each triple
changes its linking number in the course of a mirror reflection.  Thus, if
the interlacing is amphicheiral, then it must have the same 
number of
triples with linking number $+1$ as with linking number $-1$. In
particular, the total number of triples in the interlacing must be even. This
simple argument leads us to the following unexpected result.
\vskip.1in

{\bfit Theorem 1.\/} {\it If $p\equiv 3\mod4$, then every 
interlacing of $p$ lines is nonamphicheiral.}
\vskip.1in

\begin{proof} The number of triples in an interlacing of $p$ lines 
is equal to 
$p(p-1)(p-2)/6$, and this is odd if and only if $p\equiv 3
\mod4.$ 
\end{proof}
\par

\begin{figure}[htb]
\centerline{\includegraphics{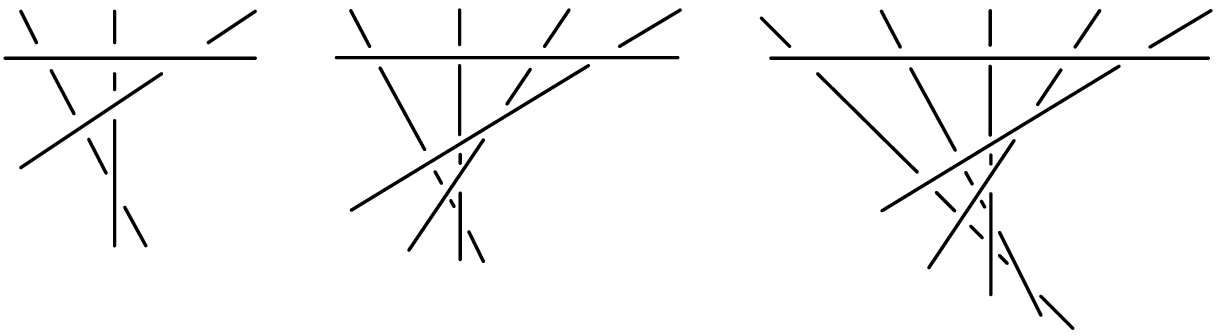}}
\caption{}
\label{f16}
\end{figure}

Theorem 1 gives an affirmative answer to the first of the four questions above.
The second question has a negative answer: for any $p\ge3$ one can construct a
nonamphicheiral interlacing of $p$ lines. This also answers question 3). The simplest
nonamphicheiral interlacings are shown in Figure \ref{f16} for $p=4,5$, 
and 6. It is easy
to continue with this sequence of examples. All of the triples of lines in the
interlacings in this sequence have the same linking number, and for this
reason we know that the interlacings are nonamphicheiral.

\par

\begin{figure}[hbt]
\centerline{\includegraphics{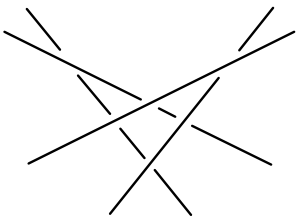}}
\caption{An amphicheiral interlacing of 4 lines.}
\label{f17}
\end{figure}

\begin{figure}[htb]
\centerline{\includegraphics{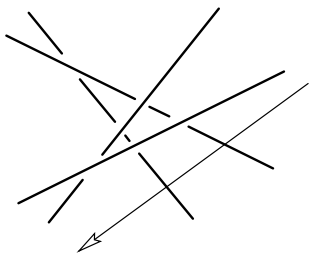}}
\caption{The same interlacing after sliding the nearest two lines to 
the left and down.} \label{f18}
\end{figure}

It remains to answer Question 4). We do not yet know whether or not 
there are
amphicheiral interlacings of $p$ lines when $p\not\equiv 3\mod4$. It is
convenient to consider separately the two cases: $p$ even, and $p\equiv1
\mod4$, although in both cases the question turns out to have a
positive answer. In Figure \ref{f17} (in which $p=4$) we show the 
simplest example of
an amphicheiral interlacing of $p$ lines with $p$ even. For any even 
number $p$, we
take two sets of $p/2$ lines, one behind the other. The lines of the set
nearest us are taken from the sequence of nonamphicheiral interlacings constructed
above (see Figure \ref{f16}). The other set of $p/2$ lines is obtained 
from the first
by rotating and then reflecting in a mirror. How do we see that the interlacing
in Figure \ref{f17} is amphicheiral? We move the set that 
is nearest us in
such a way that the part of its projection which contains all of the
intersections (in the projection) passes over and above the projection of the
other set (Figure \ref{f18}). If we then rotate Figure \ref{f18} by 
$90^\circ$ clockwise, we obtain the mirror image of the original interlacing.

\par
\begin{figure}[htb]
\centerline{\includegraphics{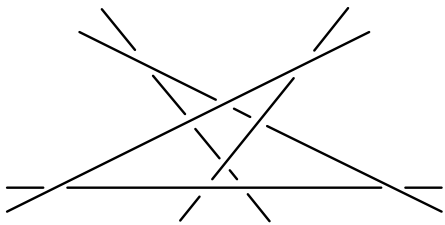}}
\caption{An amphicheiral interlacing of 5 lines.}
\label{f19}
\end{figure}

We now turn to the case $p\equiv 1 \mod4$, i.e., $p=4k+1$. 
An amphicheiral
interlacing with $k=1$ is shown in Figure \ref{f19}. Four of the lines 
form two pairs
which are situated as in the amphicheiral interlacing of four lines constructed
above. The fifth line is placed so as to separate the two lines in 
each pair. 

\begin{figure}[hbt]
\centerline{\includegraphics{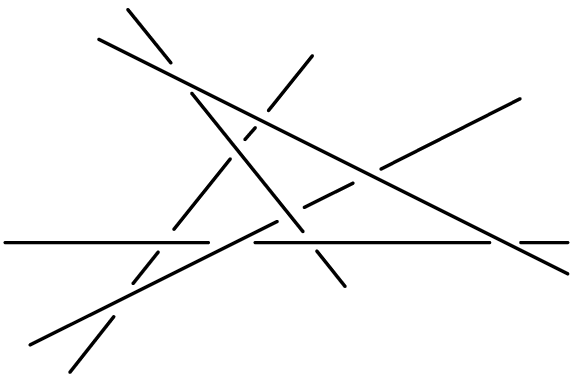}}
\caption{}
\label{f20}
\end{figure}

An isotopy between this interlacing and its mirror image can be constructed as
follows. We rotate the lines of the pair nearest us around the fifth line by
almost $180^\circ$---until the lines of the other pair are in the way (Figure
\ref{f20}). 

\begin{figure}[hbt]
\centerline{\includegraphics{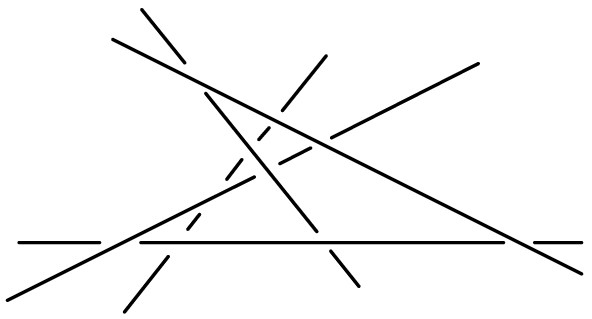}}
\caption{}
\label{f21}
\end{figure}

We then move the fifth line so that its projection passes 
to the other
side of the intersection (in the projection) of the lines that we moved before
(Figure \ref{f21}). 

\begin{figure}
\centerline{\includegraphics{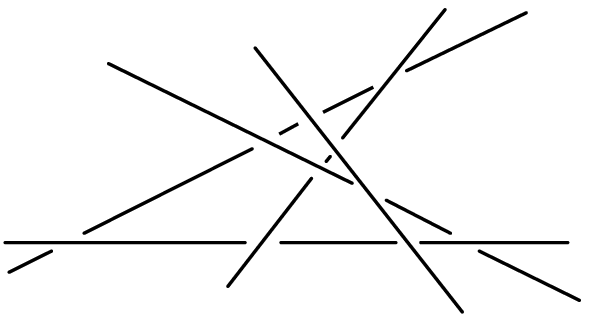}}
\caption{}
\label{f22}
\end{figure}

It remains simply to look at the resulting 
interlacing from the
opposite side. We do this by rotating it by $180^\circ$ around a vertical line
(Figure \ref{f22}). Now we see that we have the mirror image of the 
original interlacing.

\par
Using this example, it is easy to manufacture amphicheiral interlacings of $4k+1$
lines for $k>1$. Each line in Figure \ref{f19} except for the fifth 
is replaced by
an interlacing of $k$ lines which either is taken from the sequence in Figure
\ref{f16} or else is the mirror reflection of an interlacing in that 
sequence. This
must be done in such a way that the interlacings which replace the lines of one
of the pairs form an interlacing of the same type. There is no work needed to prove that the final result is an 
amphicheiral
interlacing, since the required isotopy can be obtained in the obvious way from
the one in the previous paragraph.

\section*{Four Lines}

\begin{figure}[htb]
\centerline{\includegraphics{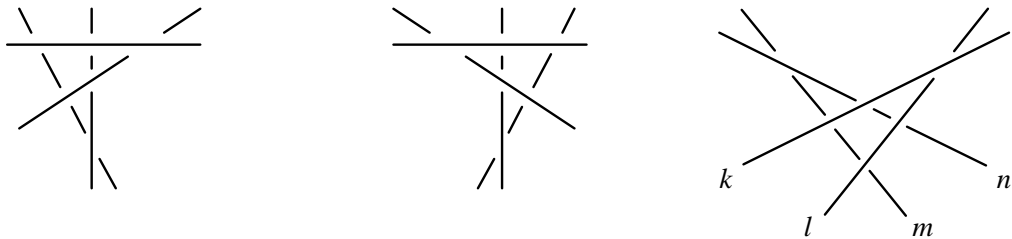}}
\caption{Interlacings of 4 lines.}
\label{f24}
\end{figure}

At this point we have actually already encountered all of the types of
interlacings of four lines. There are three of them, and they are depicted in
Figure \ref{f24}. The interlacing in Figure \ref{f16} is on the left, 
its mirror image
is in the center, and the interlacing in Figure \ref{f17} is on the 
right. We have
already proved that these three sets are not isotopic to one another: the
first one is not amphicheiral, and so is not isotopic to 
the second
one, and the third one is amphicheiral, and so is not 
isotopic to either the first or the second.
\par
In order to show that any interlacing of four lines is isotopic to one of the
interlacings in Figure \ref{f24}, we shall have to make use of the 
second of the two
objects which, as mentioned above, are associated to a triple of lines. This is
a one-sheeted hyperboloid---a surface which is usually studied in analytic
geometry. There one learns that a one-sheeted hyperboloid (henceforth
referred to simply as a hyperboloid) is made up of lines---its generatrices.
Any two generatrices in the same family are skew, while any two generatrices in
different families are either parallel or intersect. We list some other
properties of hyperboloids which we shall need:
\par
(1) if a line has three points in common with a hyperboloid, then it is a
generatrix;
\par
(2) a plane containing a generatrix of a hyperboloid intersects the hyperboloid
in two generatrices;
\par
(3) there is a hyperboloid passing through any three pairwise skew lines which
do not lie in parallel planes.
\par
These properties are simple consequences of the fact that a hyperboloid is a
surface of degree two. Of course, one could describe all of this without
appealing to analytic geometry, using the same language as the ancient Greeks,
but we shall not try the reader's patience by proceeding in that way.
\par
Thus, in order to complete the isotopy classification of four-tuples of lines,
we shall prove that any interlacing of four lines is isotopic to one of the
interlacings in Figure \ref{f24}. We take an arbitrary interlacing of 
four lines. By
moving it slightly, if necessary, we can obtain a situation where three of the
four lines (it makes no difference which three) do not lie in parallel planes.
We construct a hyperboloid through these three lines, and we observe how the
fourth line is situated relative to the hyperboloid. There are four
possibilities:
\par
(a) the line does not intersect the hyperboloid;
\par
(b) the line intersects the hyperboloid in a single point;
\par
(c) the line intersects the hyperboloid in two points;
\par
(d) the line lies on the hyperboloid.

\par
In case (d) the interlacing of four lines consists of four generatrices of the
hyperboloid, and is obviously isotopic to the left or the center interlacing in
Figure \ref{f24}.
\par
If the fourth line does not intersect the hyperboloid, then it can be brought
in toward the hyperboloid until it is tangent to the hyperboloid, i.e., the
first case can easily be reduced to case (b).

\begin{figure}[htb]
\centerline{\includegraphics{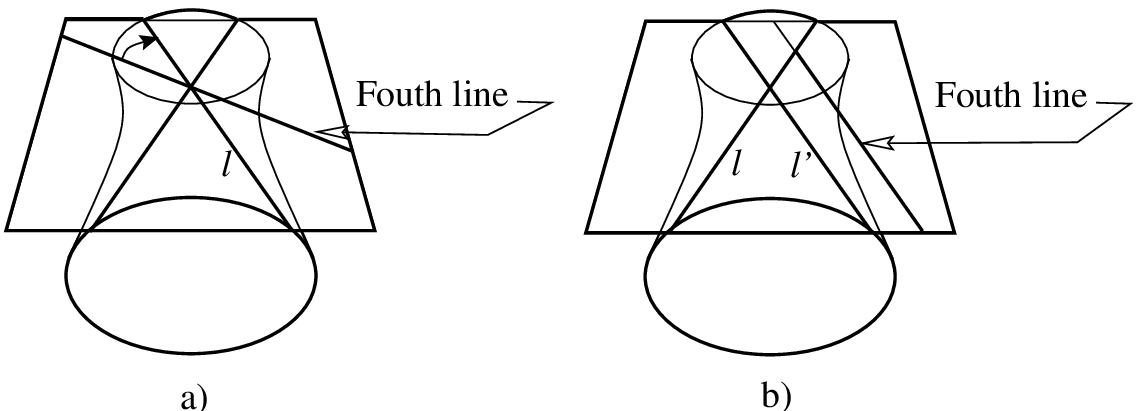}}
\caption{}
\label{f25}
\end{figure}

Case (b), in turn, reduces to either (c) or (d). To see this, we draw a
generatrix $l$ through the point of intersection of the hyperboloid with the
fourth line, where $l$ is taken in the same family of generatrices as the first
three lines of the interlacing. By property (3), the plane $a$ containing $l$
and the fourth line intersects the hyperboloid in two generatrices $l$ and
$l'$. If $l'$ intersects $l$ and the fourth line in the same point, then,
rotating the fourth line around this point of intersection in the plane $a$
until it coincides with $l$, we find ourselves in case (d) (Figure 
\ref{f25}(a)).
Otherwise, the fourth line of the interlacing is parallel to $l'$ (if this
weren't the case we would have case (c)) (see Figure \ref{f25}(b)). 
But if we perform
a small rotation toward the fourth line around the intersection point with $l$
in the plane $a$, we see that the fourth line is no longer parallel to
$l'$: it intersects $l'$, and hence it intersects the hyperboloid in two
points, giving us case (c).
\par
Now if the fourth line intersects the hyperboloid in two pints, then everything
depends on whether these points are in the same part of the hyperboloid into
which the first three lines divide it, or are in different parts (the
hyperboloid is divided into three sections). If they are in the same part, then
the fourth line can be placed on the hyperboloid without the first three lines
interfering. Then the fourth line becomes a generatrix, and we are in case (d).
If the fourth line intersects the hyperboloid in different parts, then the
interlacing is isotopic to the right interlacing in Figure \ref{f24}.

\section*{ Isotopic Lines of an Interlacing}
\par
The next step---the classification of interlacings of five lines---requires a
more careful study of the inner structure of interlacings. The reader has
undoubtedly noticed the striking difference between amphicheiral and nonamphicheiral
interlacings---compare the sets of lines in Figure \ref{f24}. The left 
and the center
interlacings both have the feature that any line of the interlacing can be
taken to any other line by means of an isotopy. This is not the case for a
amphicheiral interlacing (the right one in Figure \ref{f24}). We shall 
say that two lines of
an interlacing are {\it isotopic\/} if there exists an isotopy of the
interlacing which takes one of the lines to the other one. To be sure, strictly
speaking this is not an isotopy, because at the last moment the two lines come
together. Instead of changing the meaning of the word ``isotopy'', we are better
off leaving the meaning unchanged and adopting the following definition of
isotopic lines of an interlacing: there is an isotopy of the entire interlacing
which makes the two lines approach one another so that one can separate them
from the other lines of the interlacing by a hyperboloid (in which case there
is nothing to stop us from bringing the two lines together).
\par
Isotopic lines have the same location relative to the other lines in the
interlacing. Hence, if $a$ and $b$ are isotopic lines and $c$ and $d$ are two
other lines of the same interlacing, then the triples $a,c,d$ and $b,c,d$ have
the same linking number. Using this necessary condition for lines to be
isotopic, we can easily show that in the interlacing on the right in 
Figure \ref{f24}
the line $l$ is not isotopic to $m$. In fact, the triple $l,n,k$ has linking
number $+1$, while the triple $m,n,k$ has linking number $-1$.
\par
It is clear that, given any two isotopic lines in an interlacing, an isotopy
can be found which interchanges them and causes all of the other lines to end
up in the same place as before. Hence, isotopy of lines in an interfacing is an
equivalence relation, and the set of all lines in an interfacing is partitioned
into isotopy equivalence classes. The left and center interlacings in 
Figure \ref{f24}
each has only one equivalence class, while the right interlacing has two: the
lines $k$ and $l$ are in one class, and $m$ and $n$ are in another.
\par
If we choose one line from each equivalence class in an interlacing, then the
isotopy type of the resulting interlacing does not depend on our choice of our
choice of line in each equivalence class. This interlacing is called the {\it
derived interlacing\/}.
\par
It is useful to pass to the derived interlacing if it contains fewer lines than
the original interlacing. In order to recover the original interlacing from the
derived one, one needs a relatively small amount of additional information,
namely, how many lines were in each class and how they were linked to one
another. In fact, by means of an isotopy one can reduce the original
interlacing to a state in which the lines of each equivalence class are
generatrices of the same family on a one-sheeted hyperboloid, and the
hyperboloids
 containing the lines of the different equivalence classes do not intersect. We
leave it as an exercise to construct an isotopy that does this.
\par
The derived interlacing determines the relative location of the hyperboloids.
To recover the original interlacing it remains only to specify one of the two
families of generatrices on each hyperboloid. Here one does not have to do this
at all if the class has only one line or if there is only one class in all and
it has two lines. Otherwise the choice of a family of generatrices can be
specified by means of a numerical invariant $\vare=\pm1$ for each isotopy class
of lines in the interlacing; this is defined to be the linking number of
the triple of lines $a,b,x$, whee $a$ and $b$ are lines in the equivalence
class and $x$ is any line distinct from $a$ and $b$. We shall prove that {\it
this invariant depends only on the class of lines isotopic to $a$ and $b$\/}.
The proof will use certain formulas in which we will use the following
notation: the linking number of lines $a,b,c$ will be denoted by
$lk(a,b,c)$.
\vskip.1in
{\bfit Lemma.\/} {\it For any lines $a,b,c,d$ one has
$$lk(a,b,c)lk(a,b,d)lk(a,c,d)lk(b,c,d)=1.$$}
\vskip.1in
\par
This identity follows immediately from the definition of the linking
number of a triple of lines as the product of the linking numbers of
the three pairs of lines in the triple furnished with certain orientations. If
we give orientations to the lines $a,b,c,d$ and then compute the left side of
the above equality, we obtain the product of the squares of the linking
numbers of all possible pairs of lines $a,b,c,d.$\qed 
\par
We now prove that $lk(a,b,x)$ does not depend on $x$ when $a$ and $b$ are
isotopic lines of the interlacing. Let $y$ be any line of the interlacing which
is distinct from $a,b,x$. By the lemma we have
$$lk(a,b,x)=lk(a,b,y)lk(a,x,y)lk(b,x,y).$$
Since the  lines $a$ and $b$ are isotopic, we have $lk(a,x,y)=lk(b,x,y)$, and
hence $lk(a,b,x)=lk(a,b,y)$.
\par
It remains to show that $lk(a,b,x)$ does not depend on the choice of
representatives $a$ and $b$ of the isotopy class of lines. In fact, if $c$ is a
line which is isotopic to $a$ and distinct from $b$, then, as already proved,
we have
$$lk(a,b,x)=lk(a,b,c)=lk(a,c,b)=lk(a,c,x).\qed $$
\par
A class of isotropic lines of an interlacing whose invariant is $\vare$
 $(=\pm1)$ will be called an $\vare$-{\it class\/}.
\par
Some interlacings can be brought to the form of an interlacing of one line by
successively taking the derived interlacing. Such an interlacing is said to be
{\it completely decomposable\/}. A completely decomposable interlacing can be
characterized up to isotopy by the invariants associated with each transition
from an interlacing to its derived interlacing. We shall introduce some
notation for this characterization. An interlacing of $p$ generatrices of a
hyperboloid which form an $\vare$-class of isotopic lines will be denoted by
$\lan\vare p\ran$.
\par
We now consider $p$ hyperboloids which encompass disjoint regions and which
have the lines of the interlacing $\lan\vare p\ran$ as their axes. An
interlacing made up of $p$ subinterlacings $A_1,\dotsc,A_p$, each of which is
in the region bounded by the corresponding hyperboloid, will be denoted by 
$\lan+A_1,\dotsc,A_p\ran$ if $\vare=+1$ and $\lan-A_1,\dotsc,A_p\ran$ if
$\vare=-1$. In situations where the signs do not matter to us, we shall omit
them from the notation. For example, the interlacings in Figure 
\ref{f24} 
are characterized by the symbols $\lan+4\ran$, $\lan-4\ran$, and
$\lan\lan+2\ran,\lan-2\ran\ran$. The interlacings in Figure \ref{f16} 
are given by the
symbols $\lan+4\ran,\lan+5\ran,\lan+6\ran$. The amphicheiral interlacing of an even
number $p$ of lines that was constructed above is given by
$\lan\lan+p/2\ran,\lan-p/2\ran\ran$. In particular, the interlacing in Figure
\ref{f17} is $\lan\lan+2\ran,\lan-2\ran\ran$.
\par
Not every interlacing is completely decomposable. For example, the derived
interlacing for the interlacing in Figure \ref{f19} coincides with the 
original
interlacing, and it cannot be placed on a hyperboloid (otherwise it would not
be an amphicheiral interlacing). This is the simplest example of an 
interlacing which is not completely decomposable.

\section*{Five Lines}
\par
\begin{figure}
\centerline{\includegraphics{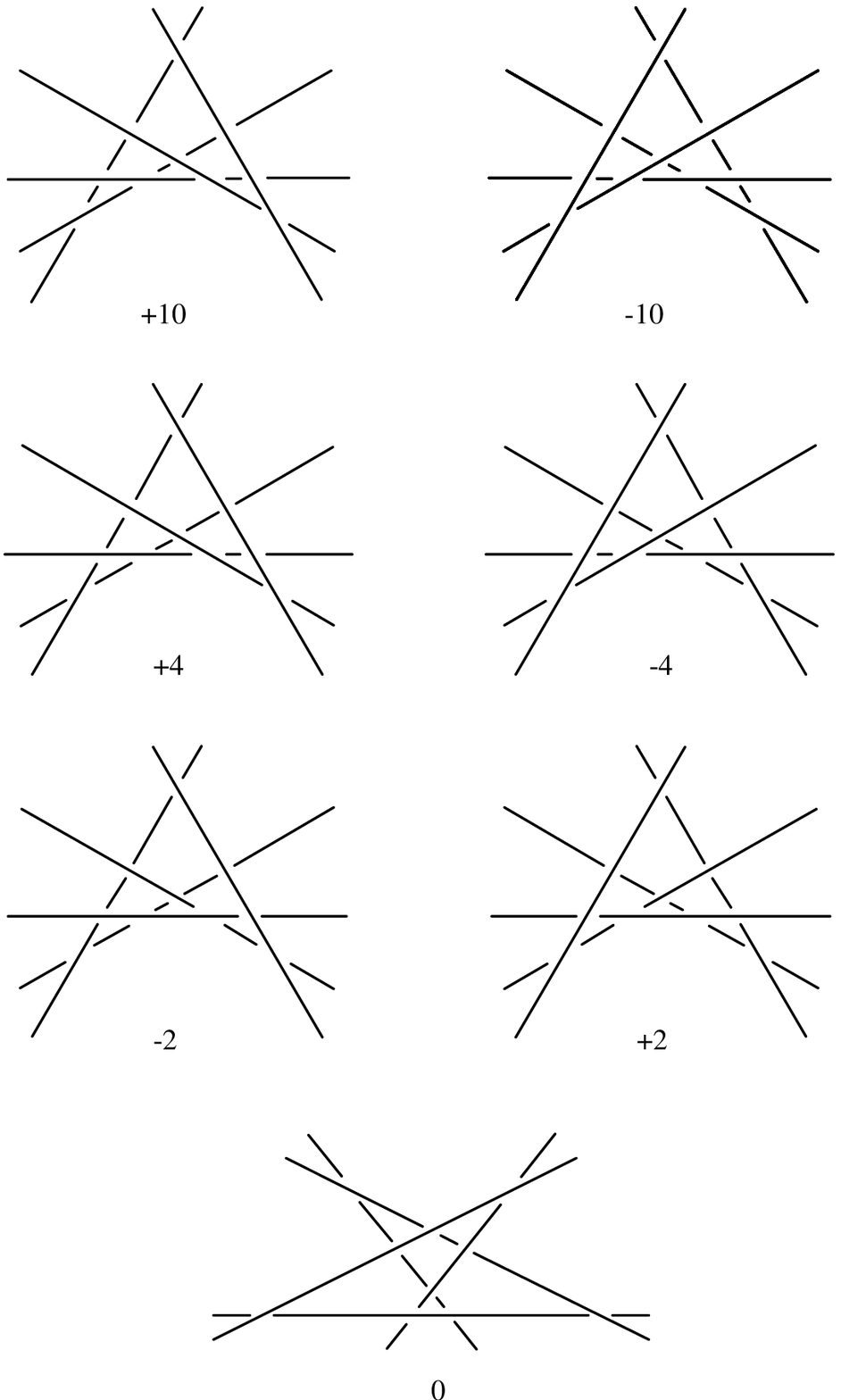}}
\caption{}
\label{f26}
\end{figure}

It can be shown (although it is not so easy as in the case of four lines) that
any interlacing of five lines is isotopic to one of the seven interlacings
shown in Figure \ref{f26}. Six of them are nonamphicheiral and 
completely decomposable; they are given by the following symbols:
\begin{align*}
&\lan+5\ran,\ \lan-5\ran,\ \lan\lan+3\ran,\lan-2\ran\ran,\
\lan\lan-3\ran,\lan+2\ran\ran,\\
&\lan+\lan1\ran,\lan-2\ran,\lan-2\ran\ran,\
\lan-\lan1\ran,\lan+2\ran,\lan+2\ran\ran.\end{align*}
The seventh is the interlacing in Figure \ref{f19}. One can prove that 
the seven
interlacings are not isotopic to one another by computing in each case the sum
of the linking numbers of the ten triples contained in the interlacing.
The results are indicated under the diagrams in Figure \ref{f26}. This 
sum is clearly
preserved under isotopy, and we see that the sums for the seven interlacings
are all different.

\section*{ Six Lines}

It is by no means so easy to show that there are in all 19 types of
interlacings of six lines (this theorem was proved by Mazurovski\u\i\ 
\cite4, \cite{5} in
1987). It is no longer possible to distinguish between nonisotopic interlacings
using only the linking numbers of the triples of lines in each
interlacing. To prove that the isotopy classes are really distinct one has to
perform computer calculations of more complicated invariants of the
interlacings. Before describing Mazurovski\u\i's basic results in more detail,
we give some definitions.
\par
We shall need a construction proposed by Mazurovski\u\i\ to characterize
interlacings of lines. Given a permutation $\si$, he constructs a corresponding
interlacing defined up to isotopy. Let $l$ and $m$ be oriented skew lines whose
linking number is $-1$.
\par
We mark off $k$ points on each line $l$ and $m$, and denote them by
$A_1,\dotsc,A_k$ and $B_1,\dotsc,B_k$, where moving form point to point with
increasing indices takes us in the direction of the line's orientation. Now,
given a permutation $\si$ of $\{1,\dotsc,k\}$, we construct an interlacing of
$k$ lines by joining $A_i$ to $B_{\si(i)}$. Following Mazurovski\u\i, we shall
denote this interlacing of the $k$ lines $A_1B_{\si(1)},\dotsc,A_kB_{\si(k)}$
by the symbol $jc(\si)$. Interlacings which are isotopic to an interlacing
constructed in this way are said to be {\it isotopy join\/}.
\vskip.1in

{\bfit Exercise} Which of the interlacings encountered above are 
isotopy joins? Show that all interlacings of five or fewer lines are 
isotopy joins.
\vskip.1in
\par
Mazurovski\u\i\ \cite6 showed that, if we want to prove that two 
interlacings
of six lines are not isotopic or if we want to determine the isotopy class of a
given interlacing of six lines, it is sufficient to use the polynomial
invariant of framed links in $\R P^3$ which was introduced by 
Drobotukhina \cite{7}. This invariant generalizes the Kauffman 
polynomial of links in
$\R^3$. 
\par
\begin{figure}[htb]
\centerline{\includegraphics{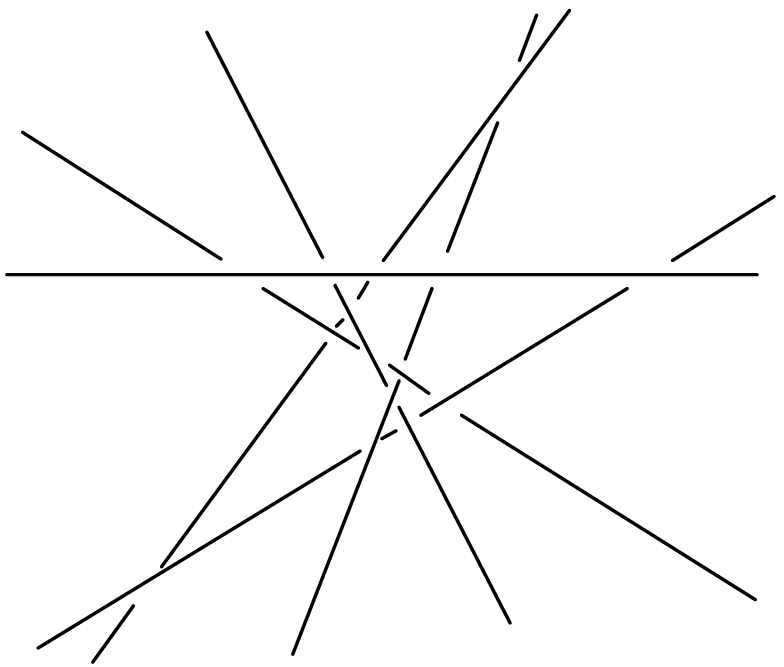}}
\caption{Interlacing $M$}
\label{f27}
\end{figure}
\begin{figure}
\centerline{\includegraphics{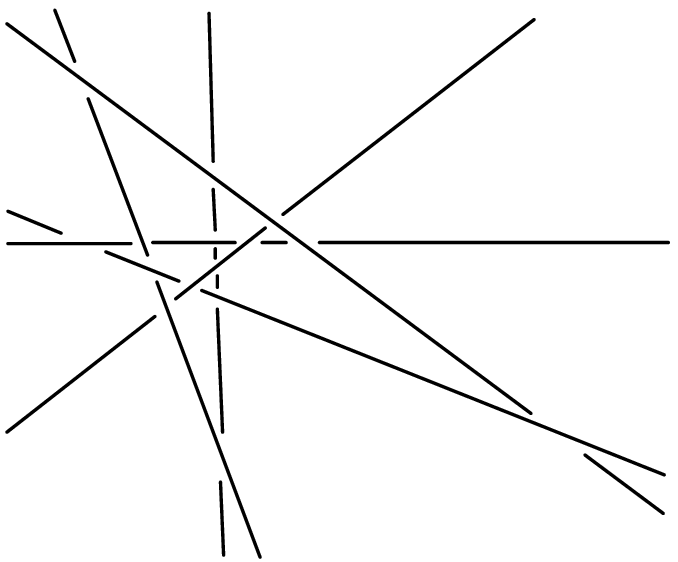}}
\caption{Interlacing $L$}
\label{f28}
\end{figure}

Return to the classification of interlacings of six lines. Of the 19 types,
15 consist of isotopy join interlacings. The remaining four are the
interlacing types $M$ and $L$ in Figures \ref{f27} and \ref{f28}, 
\noindent
and their mirror images
$M'$ and $L'$. Here $M$ and its mirror image $M'$ cannot be distinguished by
means of the linking numbers of the triples in the interlacings. But they
can be distinguished using Drobotukhina's polynomial invariant, which for $M$
is equal to
\begin{align*}
-A^{15}+6&A^{11}+6A^9-5A^7-6A^5+10A^3+16A\\
+&A^{-1}-10A^{-3}+10A^{-7}+5A^{-9},\end{align*}
and for $M'$ is equal to
\begin{align*}
5A^9+10&A^7-10A^3+A+16A^{-1}+10A^{-3}\\
-6&A^{-5}-5A^{-7}+6A^{-9}+6A^{-11}-A^{-15}.\end{align*}
Similarly, $L$ cannot be distinguished from the interlacing $jc(1,2,5,6,3,4)$
by means of the linking numbers, but these two interlacings do have
different polynomial invariants: for $L$ it is
\begin{align*}
A^{17}-5A^{13}&+15A^{9}+10A^{7}-13A^{5}-12A^{3}+15A\\
&+22A^{-1}-A^{-3}-12A^{-5}+A^{-7}+8A^{-9}+3A^{-11}\end{align*}
and for $jc(1,2,5,6,3,4)$ it is
$$A^{13}+A^{11}+4A^{7}+7A^{5}+3A^{3}+2A^{-1}+5A^{-3}+3A^{-5}
+2A^{-9}+3A^{-11}+A^{-13}.$$

The derived interlacing of $L$ coincides with $L$ itself. The same holds for
the mirror image $L'$ of $L$, the interlacings $M$ and $M'$, and also the
amphicheiral interlacing $jc(1,3,5,2,6,4)$. The interlacings $jc(1,2,4,6,3,5)$ and
$jc(5,3,6,4,2,1)$ (which are mirror images of one another) both have the same
derived interlacing, namely, an amphicheiral interlacing of five lines 
(which coincides with its own derived interlacing). The
remaining types of interlacings of six lines are completely decomposable. Of
those 12 types, two are the amphicheiral interlacings
$$\lan\lan+3\ran,\lan-3\ran\ran=jc(1,2,3,6,5,4)$$
and
$$\lan\lan-\lan1\ran,
\lan+2\ran\ran,\lan+\lan1\ran,\lan-2\ran\ran\ran=jc(1,2,4,6,5,3),$$
and the ten others can be divided into pairs of nonamphicheiral interlacings, each
pair consisting of an interlacing and its mirror image:
$$\lan-6\ran=jc(1,2,3,4,5,6),\qquad\lan+6\ran=jc(6,5,4,3,2,1);$$
\begin{align}
&\lan\lan+2\ran,\lan-4\ran\ran=jc(1,2,3,4,6,5),\\
&\lan\lan+4\ran,\lan-2\ran\ran=jc(5,6,4,3,2,1);\end{align}
\begin{align}
&\lan+\lan-3\ran,\lan-2\ran,\lan1\ran\ran=jc(1,2,3,5,6,4),\\
&\lan-\lan+3\ran,\lan+2\ran,\lan1\ran\ran=jc(4,6,5,3,2,1);\end{align}
\begin{align}
&\lan-\lan+2\ran,\lan+2\ran,\lan-2\ran\ran=jc(1,2,4,3,6,5),\\
&\lan+\lan+2\ran,\lan-2\ran,\lan-2\ran\ran=jc(5,6,3,4,2,1);\end{align}
\begin{align}
&\lan+\lan-2\ran,\lan-2\ran,\lan-2\ran\ran=jc(1,2,5,6,3,4),\\
&\lan-\lan+2\ran,\lan+2\ran,\lan+2\ran\ran=jc(4,3,6,5,2,1).\end{align}

\section*{Seven Lines}
\par
Interlacings of seven lines have been classified by Borobia and 
Mazurovski\u\i \cite{BM}. There are 74 types of interlacings of
seven lines and 48 of these types are isotopy join. As in the case of
interlacings of 6 lines, it turns out that Drobotukhina's polynomial
invariant distinguishes all the 74 types. 

A key observation which allowed Borobia and Mazurovski\u\i to obtain
this classification was a possibility to move by an isotopy each
interlacing of seven lines into a very special position. In this
position the lines are projected to a plane onto extensions of sides of
a convex polygon with seven sides, and the lines can be ordered in such
a way that the line with number $i$ passes over all the lines whose
numbers are greater than $i$.  In other words, the first line lies over
all the other lines, the second one passes over all the lines besides
the first one, the third line passes over all lines with numbers
greater than three, etc.

\section*{Interlacings of Labeled Lines}
Of course, the classification considered in the previous section, as well 
as the ones above, does not take into account any order of lines. We consider 
unordered interlacings, in which lines are not numerated or labeled. A
classification of ordered interlacings is also possible. Mazurovski\u\i
and Pavlov \cite{MP} have classified ordered interlacings of up to seven 
lines. 

When the number of lines is less than 4, the result does not
differ from the classification in unordered case. Indeed, by an
isotopy one can change arbitrarily the order of lines. 

In the case of 4 lines the number of isotopy classes of ordered
interlacings is 8. In a non-amphicheiral interlacing any two lines can
be transposed by an isotopy. Therefore the two non-amphicheiral
isotopy classes of unordered interlacings do not split when we
take into account an order of lines. So there are exactly two 
isotopy classes of ordered interlacings of 4 lines. In an amphicheiral
interlacing of 4 lines, the lines are divided into two pairs of
isotopic lines, $+1$-class and $-1$-class. Lines of the same class can
be transposed by a isotopy, while the lines of different classes
cannot. Denote the lines of the $+1$-class by $A$ and the lines
of the $-1$-class by $B$. The orderings which cannot be transformed to
each other by isotopies can be enumerated by 4-letter words made
of letters $A$ and $B$. Here are all 6 of these words: $AABB$, $BBAA$, $ABAB$,
$BABA$, $BAAB$, $ABBA$. Together with the 2 classes of non-amphicheiral
interlacings mentioned above, the corresponding 6 amphicheiral ordered
interlacings of 4 lines give totally 8 isotopy classes of ordered
interlacings of 4 lines.

Seven isotopy classes of unordered interlacings of 5 lines split
into 64 isotopy classes of ordered interlacings of 5 lines.

In the case of 6 lines, the 19 classes discussed above split into 1066.
In the case of 7 lines, 74 classes split into 43400.

\section*{Not Only Lines Can be Interlaced}
\par
We return to the definition of an interlacing of lines. We used this term to
denote a finite set of pairwise skew lines in three-dimensional space. That is,
among all possible sets of lines, we look at sets in general position which form
an everywhere dense open subset of the space of all sets of lines.
\par
The same can be done with other types of configurations. For example, we can
consider finite sets of points in three-dimensional space. We say that such a
set is {\it nonsingular\/} if for $k\le 4$ there is no set of $k$ points lying
in a $(k-2)$-dimensional subspace (i.e., a four-tuple does not lie in a plane,
a triple does not lie on a line, and all points are distinct). By an isotopy of
such a set we mean a motion in the course of which these conditions are not
violated. We say that a nonsingular set of points is {\it 
amphicheiral\/} if it is isotopic to its mirror image.
\par
We shall not treat the problem of classifying nonsingular sets of points, but
rather turn our attention to the amphicheiral problem.
\vskip.1in

{\bfit Theorem.\/} {\it A nonsingular set of $q$ points in 
three-dimensional space is
nonamphicheiral if $q\equiv 6\mod8$ or $q\equiv3\mod4$ and
$q\ge7$.}

\begin{proof} Given a nonsingular set of points, we define $s$ to be 
the sum of
the linking numbers of all triples of pairwise skew lines determined by
pairs of points in our set. If our set has $q$ points, then the number of such
triples is $\frac16\binom{q}{2}\binom{q-2}{2}\binom{q-4}{2}$. If
$q\equiv 6$ or $7\mod8$, then this number is odd, and so $s$ is also
odd, since it is a sum of an odd number of terms each of which is $\pm1$.
Clearly, $s$ is preserved under isotopies of the set of points, and it is
multiplied by $-1$ under mirror reflection. Hence, $s=0$ for an 
amphicheiral set. We
conclude that if $q\equiv 6$ or $7\mod8$, a nonsingular set of $q$
points cannot be amphicheiral. To treat the case $q\equiv 3\mod8$, 
$q\ge11$, we introduce another numerical invariant of a
nonsingular set of points. We first note that, given any two points $A$ and $B$
of our configuration, one can determine two opposite cyclic orderings of the
remaining $q-2$ points, namely, the order in which a plane rotating around the
axis $AB$ passes through them. If a triple of lines consists of the line $AB$
along with two lines joining four successive points in this ordering (more
precisely, one line joins the first point to the second and the other one joins
the third point to the fourth), then we say that the triple is {\it cyclic\/}.
Our numerical invariant of a nonsingular set of points will then be the sum of
the linking numbers of all cyclic triples of lines with distinguished
first line. If $q\ge7$, then there are $(q-2)\binom{q}{2}$ terms in this sum,
and so the sum is odd if $q\equiv 3\mod4$, $q\ge7$. On the other
hand, the sum is clearly equal to zero if we have an amphicheiral 
set. \end{proof}

\par
It is natural to ask questions about amphicheirality for nonsingular sets of
points which are analogous to the four questions discussed above in 
connection with amphicheiral interlacings of lines. We do not know
complete answers to those questions.

\par
In the same spirit one can consider a mixed situation: configurations of both
lines and points. There are various ways of defining a nonsingular 
configuration
of this type, but the most natural definition is to require that the lines in
the configuration be pairwise skew, the points not lie on the lines, and no two
points lie in a common plane with one of the lines. Even less is known about
the classification and amphicheirality of mixed configurations.
\par
When investigating problems related to geometrical objects in
Euclidean 
space,
it is often useful to extend the space to a projective space. Projective space
has even been called the ``great simplifier''. Passing to a projective space
normally enables us to find a simpler projective classification problem inside
our original classification problem, and this projective problem is usually
interesting in its own right. The case of interlacings of lines is, however, an
exception to this rule. When one goes from $\R^3$ to the projective space
$\R P^3$, an interlacing of lines corresponds to a set of disjoint
projective lines, and in this way one obtains all possible configurations of
disjoint lines in $\R P^3$ in which no line is contained in the plane at
infinity. Isotopy of interlacings is equivalent to the existence of an isotopy
between the corresponding configurations of lines in $\R P^3$ in the course
of which the lines remain disjoint.
\footnote{This is explained by the fact that, in the space of all 
configurations of
$n$ disjoint lines in $\R P^3$, the configurations containing a line in the
plane at infinity form a subset of codimension 2.}
Here one need not concern oneself with the plane at infinity. Thus, the problem
of classifying interlacings up to isotopy is actually equivalent to the
corresponding problem for configurations of lines in projective space.
\footnote{Here are two other problems which are also equivalent: the 
problem of
classifying sets of pairwise transversal two-dimensional subspaces in $\R^4$, 
with respect to motions under which they remain pairwise transversal
two-dimensional subspaces; and the problem of classifying links in the sphere
$S^3$ which are made up of great circles on the sphere, with respect to
isotopies under which the circles remain disjoint great circles on $S^3$.}
We do not get a simpler problem. But in the case of the problem of classifying
nonsingular sets of points in three-dimensional space, passing from $\R^3$
to $\R P^3$ leads to a splitting up of the problem; however, we shall not
discuss this here.

Instead we consider the following counter-part of the question on
existing of amphicheiral interlacings of a given number of lines:
For which pairs of non-negative numbers $p$, $q$ there exist an
amphicheiral nonsingular configuration of $p$ lines and $q$ points in
$\R P^3$?  We proved above the following partial results:
\begin{itemize}
\item If $q=0$, a necessary and sufficient
condition for this is $p\not\equiv3\bmod4$. 
\item If $p=0$,  it is necessary that $q\not\equiv6\bmod8$.
\item If $p=0$ and $q\ge7$, it is necessary that $q\not\equiv3\bmod4$.
\end{itemize}
The following complete answer was found by Podkorytov \cite{Pod}, 
after the previous version \cite{2} of this paper was written:
\vskip.1in

{\bfit Podkorytov Theorem.\/} (See \cite{Pod}.) {\it An amphicheiral 
nonsingular set 
of $q$ points and $p$ lines in three-dimensional real projective space
exists if and only if either 
$$q\le 3 \text{\quad  and \quad} p\equiv0\text{ or }1\bmod4, $$
or
$$q\equiv0 \text{ or }1\bmod4 \text{\quad and \quad}p\equiv0\bmod2. $$
}\qed\vskip.1in

\par
Even in the case of interlacings of lines, passing to $\R P^3$ is not
completely pointless. In $\R P^3$ we can see more clearly the topological
reasons why  interlacings are nonisotopic. As we showed at the very beginning
of the article, any interlacing can be deformed into a set of parallel lines,
and so there exists a homeomorphisms of $\R^3$ under which any interlacing
is taken to any other given interlacing with the same number of lines. In
$\R P^3$ this is no longer the case. The linking number introduced
above for oriented skew lines can be interpreted in terms of the usual linking
number in algebraic topology, applied to the corresponding lines in 
$\R P^3$ (except that we must double the topological invariant, which takes the
values $\pm1/2$, since for us the values $\pm1$ are more convenient). Moreover,
in all cases we know of, the nonisotopy of two interlacings of lines is proved
using topological invariants of the corresponding sets of projective lines in
$\R P^3$, although there probably exist nonisotopic interlacings of lines
for which the corresponding sets of projective lines can be taken to one
another by means of a homeomorphism of the ambient space.
\par
Perhaps we should show greater caution and make our definitions in accordance
with the accepted topological terminology, i.e., call interlacings of lines
isotopic if the corresponding sets of projective lines can be taken into one
another by a homeomorphism of $\R P^3$ which is isotopic to the identity
(recall that an isotopy of the homeomorphism $h\colon X\to Y$ is a family of
homeomorphisms $h_t\colon X\to Y$ with $t\in[0,1]$, $h_0=h$, such that the map
$X\times[0,1]\to Y\colon(x,t)\mapsto h_t(x)$ is continuous). Then what we
earlier called isotopies would be called {\it rigid isotopies\/}. Our cavalier
attitude about this is permissible only because at the present level of
knowledge we do not have examples of interlacings which show that these two
types of isotopies actually lead to different equivalence relations. In some
related situations, however, we do know such examples. We now discuss one such
case.

\section*{Plane Configurations of Lines}

At first glance it might seem that the world of configurations of lines in a
plane resembles the world of configurations of lines in three-dimensional
space, which we made an attempt to understand above. It is certainly easy to
give definitions for plane configurations which are analogous to the basic
definitions in this article. But, contrary to our expectations, these two worlds
have very little in common.
\par
Undoubtedly, the plane configuration analog of an interlacing of skew lines is
a configuration of lines no three of which pass through a point and no two of
which are parallel. The analog of an isotopy of interlacings is a motion during
which the lines remain lines and the conditions on the location of the lines
are preserved.
\par
Passing from the plane to the projective plane changes the problem, and here,
as usual, the projective problem turns out to be simpler and more elegant. In
the projective problem the objects are sets of projective lines in $\R P^3$
which satisfy only one condition: no three of them pass through a point. Such a
projective plane configuration of lines will be said to be {\it nonsingular\/}.
A configuration of this type can also be interpreted as a set of planes through
the origin in $\R^3$ such that no three of them contain a line.
\par
In the case of nonsingular plane configurations of lines one must distinguish
between isotopies and rigid isotopies. Two configurations are isotopic, or,
equivalently, they have the same topological type, if one can be taken to the
other by means of a homeomorphism $\R P^2\to\R P^2$. Two configurations
are said to be rigidly isotopic if they can be connected by a path in space
whose points are nonsingular plane configurations of lines.
\par
In the isotopic and rigid isotopic classification problems for plane
configurations we do not have the amphicheirality question. This is 
because the mirror
image of any configuration is isotopic to the original configuration, since a
reflection of the projective plane is isotopic to the identity map by means of
an isotopy consisting of projective transformations. (More generally, the group
of projective transformations of $\R P^2$ is connected.)
\par
The isotopic and rigid isotopic classification problems for nonsingular plane
configurations of lines have both been solved for configurations where the
number of lines is $\le7$, and in these cases the answer to both problems turns
out to be the same (see Finashin \cite8). If there are $\le5$ lines, 
the
isotopy type is determined by the number of lines. There are four types of
nonsingular plane configurations of six lines, and 11 types of nonsingular
plane configurations of seven lines. But when we reach configurations of more
than seven lines, the isotopy and rigid isotopy classifications diverge
sharply. Mnev \cite9 proved a surprising theorem, according to which, 
roughly
speaking, a set of nonsingular plane configurations of lines which are isotopic
to one another can have the homotopy type of any affine open semi-algebraic
set, and, in particular, it can have any number of connected components, i.e.,
it can contain an arbitrary number of rigid isotopy classes. (This statement is
imprecise, because in Mnev's work one considers ordered configurations in which
the first four lines are in a fixed position; otherwise one must divide out by
the action of the group of projective transformations.)
\par
The simplest example known of nonsingular plane configurations of lines which
are isotopic but not rigid isotopic can be found in Suvorov \cite{10}.
The
configurations in this example have 14 lines.

\section*{High-Dimensional Generalizations of Interlacings of Lines}

Thus, the theory of nonsingular plane configurations of lines is quite
different from the theory of interlacings of lines. This closely corresponds to
the picture one sees in the topology of manifolds: it is known that the
topology of manifolds of successive dimensions has far fewer common features
than the topology of manifolds whose dimensions differ by 4. 

In the topology of
high-dimensional manifolds one even has precise constructions which embed
various parts of $n$-dimensional topology in $(n+4)$-dimensional topology. In
surgery theory this construction is multiplication
by a complex projective plane;
in knot theory it is the two-fold covering of Bredon; and in the theory of
singularities it is the addition to our function of the sum of the squares of
two new variables. 

It seems that something similar occurs in the theory of
projective configurations. Interlacings of skew lines in three-dimensional
space appear to be related to configurations of pairwise skew
$(2k-1)$-dimensional subspaces in $(4k-1)$-dimensional space. One can define a
linking number for oriented skew $(2k-1)$-dimensional subspaces of
$(4k-1)$-dimensional space. Hence, all of the results on nonamphiheiral interlacings
that were proved using linking numbers carry over to this multidimensional
setting. 

Moreover, there is a simple construction which to any such
configuration associates a configuration of the same type with $k$ increased by
1.
This construction preserves the linking numbers, isotopic configurations
are taken to isotopic configurations, and perhaps to some extent one has an
embedding of the theory of configurations of $(2k-1)$-dimensional subspaces of
$(4k-1)$-dimensional space in the theory of configurations of
$(2k+1)$-dimensional subspaces of $(4k+3)$-dimensional space. This gives rise
to the possible development of a stable theory of projective configurations.
 
\par
Here we shall give a description of this construction. As far as we know,
it has not been published prior to the second version \cite{2} of this
paper, and it was the only original result of \cite{2}. 
Our construction of a suspension is applicable not only to configurations of
 $(2k-1)$-dimensional subspaces of $(4k-1)$-dimensional space.
It applies to any configuration  of finitely many subspaces in projective
space; it increases the subspace dimension by 2 and the ambient space dimension
by 4.
\par
We first recall the construction of the join of ordered configurations. Let
$L_1$, \dots, $L_r$ be subspaces of $\R P^p$, and let $M_1,\dotsc,M_r$ 
be subspaces of $\R P^q$. We imbed $\R P^p$ and $\R P^q$ in $\R 
P^{p+q+1}$ as skew subspaces. In the case of odd $p$ and $q$ the
imbeddings should be chosen with care about orientations: the images, 
with their native orientations, should have positive linking numbers in 
$\R P^{p+q+1}$.  Let $K_1,\dotsc,K_r$ denote the subspaces
of $\R P^{p+q+1}$ such that $K_i$ is the union of all lines which intersect
$L_i$ and $M_i$. We call the configuration of subspaces $K_1,\dotsc,K_r$ the
{\it join\/} of our two configurations.
\footnote{We have already encountered this construction. The 
isotopy join
interlacings introduced above (when we treated interlacings of six lines) are
essentially the joins of sets of points on a line.}
\par

By the {\it suspension\/} of an arbitrary configuration $L_1,\dotsc,L_r$ of
subspaces of $\R P^p$ we mean its join with a configuration of $r$
generatrices of a (one-sheeted) hyperboloid in $\R P^3$ with positive
linking number (i.e., its join with the configuration of lines in $\R 
P^3$ corresponding to the interlacing which we denoted by $\lan+r\ran)$. Since
any two lines of the interlacing $\lan+r\ran$ are isotropic, it follows that
one can find an isotopy of this interlacing which permutes the lines in an
arbitrary way. Hence, the join with an ordered configuration of subspaces
$L_1,\dotsc,L_r$ in $\R P^p$ does not depend on the order. Thus, the
suspension is well defined (up to rigid isotopy) for unordered configurations.

Two configurations of $k$-dimensional subspaces of $\R P^{2k+1}$ are
said to be {\em stably equivalent\/} if there exists $N$ such that
their $N$-fold suspensions are rigid isotopic. Mazurovskiu\i \cite{11}
has shown that this stable equivalence shares properties which are
common for various stable equivalences mentioned above. Namely, 
Mazurovski\u\i\ \cite{11} has proved that for $k>0$ any configuration 
of $\le k+5$ disjoint
$(k+2)$-dimensional subspaces of $\R P^{2k+5}$ is rigidly isotopic to 
the suspension of a configuration of
 $k$-dimensional subspaces of $\R P^{2k+1}$, and, if there are 
$\le k+2$ subspaces in the configurations, then
rigid isotopy of the suspensions is equivalent to rigid isotopy of the original
configurations of $k$-dimensional subspaces of $\R P^{2k+1}$.

This stabilization theorem was used by Mazurovski\u\i in \cite{11} for
obtaining the rigid isotopy classification of nonsingular configuration
of six $(2k-1)$-dimensional subspaces in $\R P^{4k-1}$. He proved that
when $k>1$ such configuration is defined up to rigid isotopy by the
linking numbers. Recall that this is not the case for $k=1$, that is
for configurations of lines in the 3-space. Suspension makes
configuration $M$ shown in Figure \ref{f26} and its mirror image $M'$
rigidly isotopic. Recall that $M$ and $M'$ are distinguished by the
Kauffman bracket polynomial. Thus, there is no generalization of
the Kauffman bracket to high-dimensional nonsingular configurations
which would be preserved under suspension. 

Then Khashin and Mazurovski\u\i \cite{KM} proved that 

{\em Two nonsingular
configurations of $k$-dimensional subspaces of $\R P^{2k+1}$ are stably
equivalent if and only if they have the same linking numbers of the
subspaces.} 

This means that there exists a bijection between the set of
the $k$-subspaces of the first configuration and the set of the
$k$-subspaces of the other configuration such that the linking numbers
of the corresponding subspaces are equal. 

Algebraic techniques developed for that was used in \cite{KM} also for
obtaining the following two results about interlacings of skew lines in
the 3-space: 

{\em  Two isotopy join interlacings are rigidly isotopic if and only
if they have the same linking numbers (and hence stably equivalent).}

{\em An interlacing of skew lines which has the same linking numbers as
the configuration of disjoint generatrices of a (one-sheeted) hyperboloid in 
$\R P^3$  is rigidly isotopic to this configuration of generatrices.}

\section*{Connection with Real Algebraic Surfaces of Degree 
$4$} 

Almost everything in the first two-thirds of the article concerning 
interlacings of
lines, as well as everything concerning nonsingular sets of points in
three-dimensional space, was published by Viro in 1985 in the note 
\cite3.
Interest in this subject was stimulated by work of Kharlamov on the
classification of nonsingular real projective algebraic surfaces of degree 4 up
to rigid isotopy (by which one means isotopies consisting of nonsingular
algebraic surfaces). Earlier, a coarser classification of such surfaces up to
mirror reflections and rigid isotopies was found by Nikulin 
\cite{12}; and
Kharlamov, using a very complicated technique which involved passing to the
complex domain, proved that certain surfaces are nonamphiheiral, in the sense that
they are not rigid isotopic to their mirror images. It would be worthwhile to
find an elementary proof.
\par
Some of these surfaces decompose in the ambient three-dimensional space into a
one-sheeted hyperboloid with handles and a number of separate spheres (the sum
of the number of handles and the number of spheres is at most ten, and there
are other restrictions, but we shall not dwell on this). From Harnack's theorem
on the number of components of a plane curve it follows that every plane
intersects at most three of the spheres of this surface. Hence, if we choose
one point on each sphere, we obtain a nonsingular set of points whose isotopy
type is determined by the surface, and a rigid isotopy of the surface
corresponds to an isotopy of the set of points. Thus, if there are six or seven
spheres, the surface must be nonamphiheiral. In a similar way Kharlamov proved that
many other degree 4 surfaces are nonamphiheiral and completed the classification of
nonsingular surfaces of degree 4 (see \cite{13}). However, he was able 
to prove
that certain of the surfaces are nonamphiheiral only by passing to the complex
domain and using the full theory of K3-surfaces.

\end{document}